\documentclass{amsart}
\usepackage{amssymb}
\usepackage{graphics}

\theoremstyle{plain}
\newtheorem{lemma}{Lemma}
\newtheorem{proposition}{Proposition}

\theoremstyle{definition}
\newtheorem{definition}{Definition}

\theoremstyle{definition}
\newtheorem{remark}{Remark}

\begin{document}
\title[Computation of Galois groups]
     {Computation of  Galois groups associated to  \\ the 2-class towers  of 
	some quadratic fields }
\author{M. R. Bush}
\address{Mathematics Department \\
	 University of Illinois \\
	 Urbana-Champaign, IL 61801 }
\email{mrbush@math.uiuc.edu}
\thanks{The author would like to thank Nigel Boston for suggesting
this problem and for many helpful discussions.}
%11R29 -- class groups, class numbers etc.
\begin{abstract} The $p$-group generation algorithm from computational
group theory is used to obtain information about large quotients of the 
pro-$2$ group $G = \text{Gal}\,(k^{nr,2}/k)$ for $k = \mathbb{Q}(\sqrt{d})$
with $d = -445, -1015, -1595, -2379$. In each case we are able to
narrow the identity of $G$ down to one of a finite number of  explicitly
given finite groups. From this follow several results regarding the 
corresponding 2-class tower.
\end{abstract}
\maketitle

\section{Introduction} In recent years there has been much interest
in trying to determine the behaviour of the $p$-class tower of a
quadratic field especially in the case where $p = 2$. See for instance
 \cite{BLS1}, \cite{BLS2}, \cite{BG1}, \cite{Lem1} and \cite{Yam}.
Note that by {\em ($p$-)class tower of $k$} we mean the chain of fields
\[  k = k_0 \subseteq k_1 \subseteq \ldots \subseteq k_n \subseteq \ldots \]
where $k_{n+1}$ is the Hilbert ($p$-)class field of $k_n$ for each nonnegative
integer $n$. We say that the tower is {\em finite} if $k_n = k_{n+1}$ for
some $n$, and infinite otherwise. If it is finite then the minimal such $n$
is called the {\em length} of the tower.
In this paper we  show that the $2$-class towers of four
imaginary quadratic fields $k$ are finite 
using a computational method first introduced by Boston and 
Leedham-Green \cite{BG1}. In fact we are able to
give a short list of candidates for the Galois group
$\text{Gal}\,(k^{nr,2}/k)$ in each case. Here $k^{nr,2}$ denotes the 
maximal unramified $2$-extension of $k$ which is obtained by taking the union
of the fields occurring in the 2-class tower of $k$.

The first field we consider is $k = \mathbb{Q}(\sqrt{-2379})$.
It was noted by Stark \cite{St} that this field has root discriminant
$\sqrt{|d_k|} \approx 48.8$ which is just above the best known lower
bound ($\approx 44.7$ under GRH) for $R = \liminf_{n \rightarrow \infty} R_n$ where
$R_n$ is the minimal root discriminant over all imaginary number fields of degree $n$
(see \cite{Odl}). If the class tower above $k$ were infinite we would thus
obtain a fairly tight upper bound on possible sharpenings of this lower bound
(since root discriminants remain constant in class towers). Currently the best 
known upper bound on this lower bound
is around $83.9$ (see \cite{HM1}) so this would be a significant improvement.
We do not resolve whether or not the class tower of $k$ is infinite here
but our methods do show that the $2$-class tower is finite of length 2. 

The three other fields we consider are $k = \mathbb{Q}(\sqrt{-445})$
, $\mathbb{Q}(\sqrt{-1015})$ and $\mathbb{Q}(\sqrt{-1595})$. 
It is observed in \cite{BLS1} that under GRH these fields
are the first examples of imaginary fields with finite 
(2-)class towers and $\text{rank}\,\text{Cl}_2(k_1) \geq 3$.
Here we determine that each field has finite $2$-class tower unconditionally.
We are also able to obtain the exact length of the tower and the 
2-class groups which occur in it in each case.

\section{The Method}
The method is based on the $p$-group generation algorithm introduced 
by O`Brien \cite{Obrien}. We now recall some definitions. Let $G$ be a
pro-$p$ group. We recursively define a series of closed subgroups of
$G$
\[  G = P_0(G) \geq P_1(G) \geq P_2(G) \geq \ldots \]
by setting $P_n(G) = P_{n-1}(G)^p [G,P_{n-1}(G)]$ for each $n \geq 1$.
Here the group on the righthand side is the closed subgroup generated
by all $p$ th powers of elements in $P_{n-1}(G)$, and commutators of
elements from $G$ and $P_{n-1}(G)$. If $G$ is a finite $p$-group then the series above
is finite and the smallest $c$ such that $P_c(G) = \{1\}$ will be called
the {\em $p$-class} of $G$. A $p$-group $H$ is called a {\em descendant} of
$G$ if $H / P_c(H) \cong G$ where $c$ is the $p$-class of $G$. It is an
{\em immediate descendant} if it has $p$-class $c + 1$. The $p$-group
generation algorithm finds representatives (up to isomorphism) of 
all the immediate descendants of a given finite
$p$-group $G$. 

Now, fix an ordered pair $(G,\{G_i\}_{i=1}^{n})$ where 
$G = \text{Gal}\,(k^{nr,2}/k)$ is a pro-$2$ group
and $G_i$ is a closed subgroup of $G$ for $i = 1,\ldots\,,n$. We will be
interested in the pairs $G_{(m)} = (G/P_m(G), \{\overline{G_i}\}_{i=1}^{n})$
 where $m \geq 0$ and $\overline{G_i}$ denotes
the image of the subgroup $G_i$ under the natural map to the quotient 
$G/P_m(G)$. We note that $G/P_m(G)$ is always finite.
\begin{definition} A pair
$(H,\{H_i\}_{i=1}^n)$ will be called a {\em representative} of the pair
$G_{(m)}$ if there exists an isomorphism
$\psi : H \longrightarrow G/P_m(G)$ such that $\psi(H_i) = \overline{G_i}$
for each $i = 1,\ldots\,,n$. 
\end{definition}
%Observe that any two representatives
%$(A,\{A_i\}_{i=1}^n)$ and $(B,\{B_i\}_{i=1}^n)$ of
%the same pair are {\em equivalent} by which we mean there exists an
%isomorphism $\phi:A \longrightarrow B$ such that $\phi(A_i) = B_i$ 
%for all $i = 1,\ldots\,,n$. 
We now have a lemma which follows easily from the definitions.

%Our aim here is to give a method for computing a list
%of pairs containing a representative of the pair 
%$G_{(m)} = (G/P_m(G), \{\overline{G_i}\}_{i=1}^{n})$ for each $m > t$ when we
%have been given a list of pairs containing a representative for $m = t$.

\begin{lemma}\label{lem1}
Suppose that $(P,\{P_i\}_{i=1}^n)$ is a representative of $G_{(m)}$ and
 $(Q,\{Q_i\}_{i=1}^n)$ is a representative of $G_{(m+1)}$, then
\begin{enumerate}
\item $Q$ is an immediate descendant of $P$.
\item There exists a surjective map $f:Q \rightarrow P$ such that
  $f(Q_i) = P_i$ for $i = 1,\ldots\,,n$.
\end{enumerate}
\end{lemma}

Given a representative of $G_{(m)}$ Lemma \ref{lem1} allows us to compute a
finite list of pairs containing a representative of $G_{(m+1)}$.
This is because a finite group $P$ has only finitely many immediate 
descendants $Q$ up to isomorphism and for  each  descendant $Q$ there are
only finitely many surjective maps $f:Q \rightarrow P$. Thus,
starting with a list of pairs known to contain a representative of
$G_{(t)}$ for some $t$ we can (by applying Lemma \ref{lem1} repeatedly)
compute a list of pairs which must contain a representative of 
$G_{(m)}$ for any $m > t$.

As it stands the above is not particularly useful since the
number of pairs grows very quickly with $m$. To try
to eliminate some of these we now suppose that  
we know the direct product decomposition of the abelian group $G_i/G_i'$ into
cyclic groups for $i = 1,\ldots\,, n$. From now on this information about the
decomposition  will be referred to as the {\em abelian quotient invariants of
$G_i$}. It will usually be given by listing the orders of the cyclic subgroups
involved, so for instance if $G_i/G_i' \cong C_2 \times C_4 \times C_4$ (where
$C_2$ and $C_4$ are cyclic groups of orders 2 and 4 respectively) then
we would say that $G_i$ has abelian quotient invariants $[2,4,4]$.
We note that in practice it is possible to get hold of such information
for low index subgroups of $G$ by computing the 2-Sylow subgroups of the 
ideal class groups of small subextensions of $k^{nr,2}/k$ and then 
applying Class Field Theory.
If we also assume that $G_i \supseteq P_t(G)$ for each $i$ then
for $m \geq t$ the abelian quotient invariants of the image of $G_i$
in $G/P_m(G)$ must be a quotient of those of $G_i$.
In the examples considered below the application of 
these additional restrictions
drastically reduces the number of candidates for representatives
of $G_{(m)}$ for each $m$ allowing us to obtain useful information
about the extension $k^{nr,2}/k$.

In summary, given a list of pairs $\mathbf{L}_m$ containing a representative of 
$G_{(m)}$ we compute a list $\mathbf{L}_{m+1}$containing a representative 
of $G_{(m+1)}$ as follows
\begin{enumerate}
\item For each pair $(P,\{P_i\}_{i=1}^n)$ in $\mathbf{L}_m$ we compute a list
of all the immediate descendants of $P$ (up to isomorphism) using the
$p$-group generation algorithm.
\item For each immediate descendant $Q$ we construct all possible
surjective maps $f:Q \rightarrow P$.
\item For each $Q$ and $f:Q \rightarrow P$ we check to see whether the 
abelian quotient invariants of $f^{-1}(P_i)$ are a quotient of those of
$G_i$ for every $i = 1,\ldots\,,n$. If this is the case then the
pair $(Q,\{f^{-1}(P_i)\}_{i=1}^n)$ is appended to the list $\mathbf{L}_{m+1}$.
\end{enumerate}
In practice we can usually compute a suitable starting list $\mathbf{L}_1$ or $\mathbf{L}_2$. In the examples
considered in the next section the sequence of lists $\mathbf{L}_m$ that are then generated 
terminate, ie. $\mathbf{L}_m$ is empty for sufficiently large $m$.
This implies  that $G = \text{Gal}\,(k^{nr,2}/k)$  actually occurs as a group 
in one of the pairs on these lists and so must be finite.
We now note that if $(H,\{H_i\}_{i=1}^n)$ is  a 
representative of $(G,\{G_i\}_{i=1}^n)$ then 
the abelian quotient invariants of $H_i$ must match  those of $G_i$ exactly.
This observation motivates the following definition.
\begin{definition}
A finite group $H$ will be called a {\em candidate for $G$} if there exists a pair $(H,\{H_i\}_{i=1}^n)$  
such that the abelian quotient
invariants of $H_i$ are equal to those of $G_i$ for all $i = 1,\ldots\,,n$.
\end{definition} 
As we generate pairs we keep track of all the candidates for $G$ that we find since
in the event the sequence of lists terminate we know that that $G$ must be isomorphic
to one of them.

\section{Some Examples}
Two packages were used to carry out the actual computations in this section.
KASH (see \cite{kash}) was used  to construct $2$-class fields above $k$ up to degree
16, as well as for some ideal class group computations. MAGMA (see \cite{magma}) was
used to calculate ideal class groups, as well as determining the subfield
lattices and Galois groups of certain degree 16 extensions of $\mathbb{Q}$. We also made
use of MAGMA's implementation of the $p$-group generation algorithm.
Throughout this section $G = \text{Gal}\,(k^{nr,2}/k)$.

\subsection{The field $\mathbf{k = \mathbb{Q}(\sqrt{-2379})}$.}  
Using KASH we have $G/G' \cong [4,4]$.
% (note that $[4,4]$  stands for the direct product of two cyclic groups of order 4)
It follows from this that $G/P_1(G) \cong [2,2]$.
Calculation of some $2$-class fields shows that 
$\mathbb{Q}(\sqrt{-3},\sqrt{13},\sqrt{61})$ is a degree
4 subextension of $k^{nr,2}/k$ with Galois group over $k$ equal to
$[2,2]$. This implies that the subgroup $P_1(G)$ of $G$ corresponds
to the subfield $\mathbb{Q}(\sqrt{-3},\sqrt{13},\sqrt{61})$.
 The following table shows
the intermediate fields together with the abelian invariants of
the 2-Sylow subgroup of the corresponding ideal class group. \\[10pt]
\centerline{
\renewcommand{\arraystretch}{1.25}
\begin{tabular}{|l|l|}
\hline
$k = \mathbb{Q}(\sqrt{-3\cdot 13 \cdot 61})$ & $[4,4]$ \\ \hline
$k(\sqrt{61})$ & $[2,2,8]$ \\ \hline
$k(\sqrt{13})$ & $[2,2,8]$ \\ \hline
$k(\sqrt{-3})$ & $[2,2,16]$ \\ \hline
$\mathbb{Q}(\sqrt{-3},\sqrt{13},\sqrt{61})$ & $[4,4,8]$ \\ \hline
\end{tabular}
}
\newline

\begin{figure}[h]
 \scalebox{.6}{\includegraphics{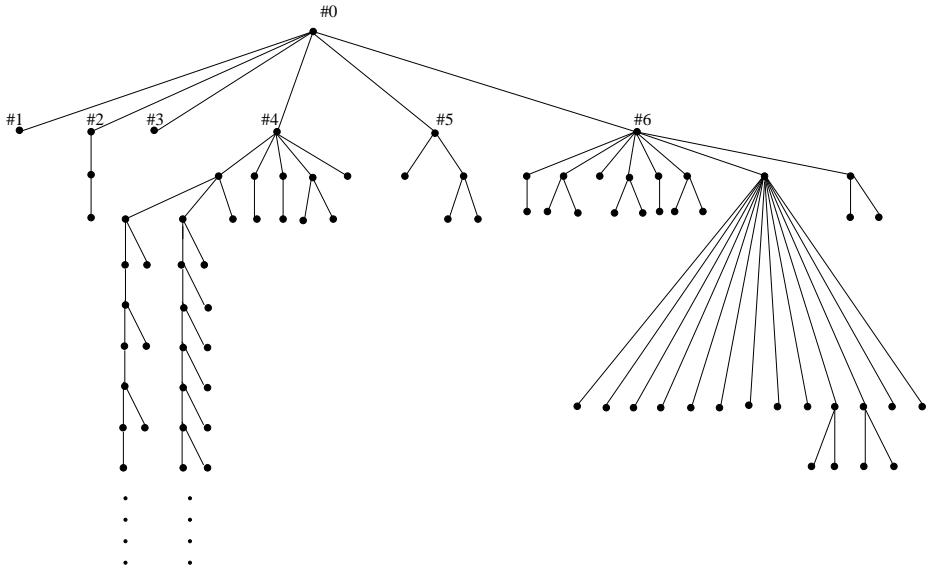}}
 \caption{A tree of groups}\label{fig1}
\end{figure}

We can thus construct a representative of $G_{(1)} =
(G/P_1(G),\{\overline{G_i}\}_{i=1}^5)$ where the groups $G_i$ 
are the subgroups of $G$ corresponding to the intermediate fields of the
extension $\mathbb{Q}(\sqrt{-3},\sqrt{13},\sqrt{61})/ k$.
Let $\mathbf{L}_1$ be the list containing just this pair.
We now apply the method from
Section 2 using the abelian quotient invariants given in the table
to generate a sequence of lists $\mathbf{L}_m$ for $m > 1$. The
collection of groups that occur in the pairs in these lists can
be represented by a tree as shown in Figure \ref{fig1}.

In this picture we have only displayed vertices which  represent groups having
 at least one immediate descendant and which occur in one of the lists
$\mathbf{L}_m$. Groups which do not have any descendants but which occur
 in one of the lists are not shown. This is done primarily for convenience since
including such (terminal) groups would clutter up the picture and we are really
only interested in determining whether or not this and other such trees are finite.
Note that in this tree a vertex $Q$ is a child of 
another vertex $P$ if the group corresponding to $Q$ is an immediate
descendent of the group corresponding to $P$. Thus the vertex numbered
$0$ corresponds to the group $[2,2]$ of $2$-class 1 (from the pair in 
$\mathbf{L}_1$). The groups numbered 1 to 6  have $2$-class 2
and come from the pairs in $\mathbf{L}_2$ and so on. Note that the groups
in each level of the tree are pairwise nonisomorphic but  
each one may occur in more than one pair on the corresponding list.

While generating the lists a total of 81 candidates for $G$ were found (all of them
descendants of the group corresponding to vertex 6). Most of these are not displayed
in Figure \ref{fig1} since they do not have any descendents. At this stage 
we cannot conclude that $G$ must be one of these groups since after
several computations part of the tree (below vertex 4) is still growing.
We have however gained some useful information about the quotient 
$G/P_2(G)$. In particular since no candidates for $G$ were found in the 
finite subtrees lying below vertices 1,2,3 and 5 we can conclude
that $G/P_2(G)$ must be isomorphic to the group corresponding to 
either vertex 4 or vertex 6. Let us denote these two groups by $H_4$ and
$H_6$ respectively. $H_4$ is generated by $\{x_i\}_{i=1}^4$ subject to
the power commutator presentation
\[ H_4\, : \, x_1^2 = x_4 \qquad  [ x_2 , x_1 ] = x_3  \]
$H_6$ is generated by $\{x_i\}_{i=1}^5$ subject to the power commutator
presentation
\[ H_6\, : \, x_1^2 = x_4 \quad  x_2^2 = x_5 \quad [ x_2 , x_1 ] = x_3  \]
\begin{remark} \label{pcrem}
Note that in any power commutator presentation if a power $x_r^2$ 
or commutator $[x_r,x_s]$ does not occur on the left hand side
of the given relations then it is assumed to be trivial.
\end{remark}

\begin{figure}[h]
 \scalebox{.55}{\includegraphics{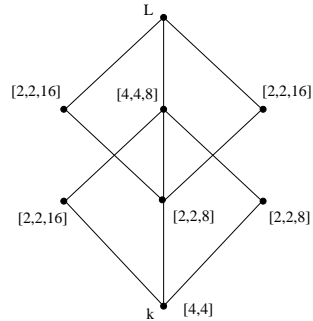}}
 \caption{A lattice of subfields and 2-class groups}\label{fig2}
\end{figure}

Using KASH we can find a degree 16  Galois extension $L$ over $\mathbb{Q}$ 
unramified over $k$. One such field is defined by the polynomial 
\begin{multline*}
 x^{16} - 2158 x^{14} - 1166 x^{13} + 1886402 x^{12} + 1125558 x^{11} - 738996514 x^{10} 
 \\+  24633036 x^9   
    + 88589769625 x^8 - 114401828130 x^7 + 12435312336118 x^6 
   \\ + 15732271973132 x^5  
+ 506694031967064 x^4 -    98005626098698 x^3 \\ + 
10557300816504844 x^2 
- 7195589177918350 x + 41648817878658175
\end{multline*}
and its Galois group is the direct product of the dihedral group
of order 8 and the cyclic group of order 2.
By considering the lattice of subfields of this extension over $k$
it can be shown that $\text{Gal}\,(L/k) \cong [2,4]$. 
The lattice of subfields and the 2-Sylow subgroups of the corresponding
ideal class groups are shown in Figure \ref{fig2}.

Now $[2,4]$ is a group of $2$-class 2 which means that the subgroup $P_2(G)$ 
of $G$ corresponds to some intermediate field of the extension $k^{nr,2}/L$. 
So we consider the
pair $G_{(2)} = (G/P_2(G), \{\overline{G_i}\}_{i=1}^7)$ where $G_i$ are the subgroups of $G$
corresponding to the intermediate fields (of degree at most 4) in the lattice above.
A list $\mathbf{L}_2$ containing a representative of $G_{(2)}$ can be constructed by
taking each group $H_r$ ($r = 4$ or $r = 6$) and adjoining all pairs of the form
$(H_r, \{f^{-1}(R_i)\}_{i=1}^7)$ where $f$ runs through all the surjective homomorphisms
from $H_r$ onto $[2,4]$ and $\{R_i\}_{i=1}^7$ are the subgroups of $[2,4]$ of index at most
4. When the method from Section 2 is applied to $\mathbf{L}_2$ 
the sequence of lists generated terminates.
The extra information about the abelian quotient invariants of the two additional subgroups
of index 4 is restrictive enough to show $G$ must be finite. The number of candidates for $G$ 
found (with the additional 
restrictions) also decreases from 81 to 24. An investigation of these 24 groups shows that
they fall naturally into two classes based on the abelian quotient invariants of their
subgroups of index 4. In one of the classes (consisting of 8 groups) each group has
several subgroups of index 4 with abelian quotient invariants $[4,32]$. None of the 
subgroups of index 4 of the groups in the other class have these invariants. KASH
 can be used to show that the extension of $\mathbb{Q}$ defined by the polynomial
\[ x^8 + 9494 x^6 + 33992937 x^4 + 54094064336 x^2 + 32175758727424   \]
%list[16]  Finished/C39_61/kashdata. 
%\begin{multline*}
%x^{16} - 4295 x^{14} + 1097560705 x^{12} - 356243647970 x^{10} + 27711347664766 x^8 \\
%+  198617618065150 x^6 + 431840387310625 x^4 + 7872666015625 x^2 + 152587890625 
%\end{multline*}
%is an unramified extension of $k$ with Galois group of 2-class 2 and with a degree 4 
%subextension over $k$ whose ideal class group is $[4,32]$. 
is an unramified extension of $k$ with ideal class group $[4,32]$.
We conclude that
$G$ must be isomorphic to one of 8 possible groups. Each one of these can be obtained
by selecting $r,s,t \in \{ 0,1 \}$ and then considering the group generated by
$\{x_i\}_{i=1}^{11}$ subject to the power commutator presentation 
\begin{align*}
     &[x_2,x_1] = x_3 \hspace*{2.5cm}  &[x_4,x_3] &= x_{10}\hspace*{2.5cm} &[x_6,x_1] &=  x_9\\
     &[x_3,x_1] = x_6  &[x_5,x_1] &= x_6 x_7 x_8 x_9 x_{10} &[x_8,x_1] &=  x_{10} \\ 
     &[x_3,x_2] = x_7  &[x_5,x_3] &= x_{10} x_{11} &[x_8,x_2]  &=  x_{10} x_{11}  \\
     &[x_4,x_2] = x_8  &[x_5,x_4] &=  x_{10} x_{11} &[x_9,x_1] &=  x_{11} \\
      &x_1^2 = x_4 &x_5^2 &= x_6  x_9  x_{10}^s x_{11}^t  & &\\
      &x_2^2 = x_5 &x_6^2 &= x_9  x_{10}  x_{11} & &\\
      &x_3^2 = x_6  x_8  x_9  x_{10}  &x_7^2 &= x_{10}  x_{11} & &\\
      &x_4^2 = x_7 x_{11}^r    &x_9^2 &= x_{11} & &  
\end{align*}
%   old::   \quad \makebox[2cm][l]{$\text{where}\: s,t \in \{0,1\}$}
Remark \ref{pcrem} about power commutator presentations applies here also.

The 8 groups defined by these presentations are very similar. They all have order
$2^{11}$ and 2-class $5$. They also possess derived series of the same length and with 
the same abelian factors. This allows us to deduce 
\begin{proposition} The $2$-class tower of the imaginary quadratic field
$k = \mathbb{Q}(\sqrt{-2379})$ is finite of length
2, ie. $k = k_0 \subset k_1 \subset k_2 = k^{nr,2}$. We have
$\text{Gal}\,(k_1/k_0) \cong [4,4]$ and $\text{Gal}\,(k_2/k_1) \cong [2,4,16]$.
\end{proposition}

\subsection{The field $\mathbf{k = \mathbb{Q}(\sqrt{-445})}$.}  
In this case $G/G' \cong [2,4]$ so we still have $G/P_1(G) \cong [2,2]$.
The field $\mathbb{Q}(\sqrt{-1},\sqrt{5},\sqrt{89})$ is a degree 4 subextension
of $k^{nr,2}/k$. We can calculate the 2-Sylow subgroup of 
the ideal class group of each intermediate field, and construct a representative
of the pair $G_{(1)} = (G/P_1(G),\{\overline{G_i}\}_{i=1}^5)$ where the groups $G_i$ 
are the subgroups of $G$ corresponding to the intermediate fields.
% of the extension $\mathbb{Q}(\sqrt{-1},\sqrt{5},\sqrt{89})/\mathbb{Q}$.
Applying the method from Section 2 we are able to narrow down $G/P_2(G)$ to
one of 3 groups. Two of these can be identified as $[2,4]$ and $D_4$ (the 
dihedral group of order 8). These can be eliminated as possibilities since
it is easy to find examples of degree 8 subextensions of $k^{nr,2}/k$ with either of
these groups as the Galois group. This implies that both groups must arise
as quotients of $G/P_2(G)$. It follows that $G/P_2(G)$ must be isomorphic to
the remaining group of order 16 which we will denote $H_3$. The group 
$H_3$ is actually isomorphic to the group $H_4$ in the previous example.

\begin{figure}[h]
 \scalebox{.55}{\includegraphics{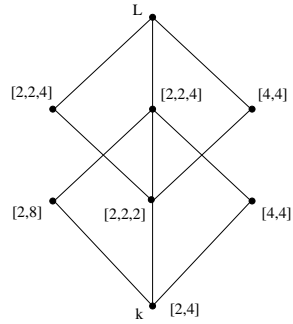}}
 \caption{A lattice of subfields and 2-class groups}\label{fig3}
\end{figure}

As before we now consider an unramified degree 8 extension $L$ over $k$ defined over
$\mathbb{Q}$ by the polynomial
\begin{multline*}
 x^{16} + 12 x^{14} + 4554 x^{12} + 17928 x^{10} + 
    2231251 x^8  + 13625880 x^6 \\ - 10866150 x^4 - 143437500 x^2 + 244140625
\end{multline*}
We have $\text{Gal}\,(L/k) \cong[2,4]$. The lattice of subfields of $L/k$ 
together with the 2-Sylow subgroups of the corresponding
ideal class groups can be seen in Figure \ref{fig3}.
We construct $\mathbf{L}_2$ by adjoining all pairs of the form
$(H_3,\{f^{-1}(R_i)\}_{i=1}^7)$ where $f$ runs over all the surjective maps
from $H_3$ to $[2,4]$ and $\{(R_i)\}_{i=1}^7$ are the subgroups of
$[2,4]$ of index at most $4$. When we apply the method starting with
$\mathbf{L}_2$ the sequence of lists generated terminates and a total
of 12 candidates for $G$ are found. To try and determine which one is isomorphic to
$G$ we look again at the abelian quotient invariants of their index 4
subgroups. We find that only 2 of the 12 groups have an index 4 subgroup with
invariants $[2,16]$. KASH can be used to show that the extension of $\mathbb{Q}$
defined by the polynomial 
\[ x^8 + 702 x^4 + 130321  \] 
%list[16]  C5_89ck/kashdata
%\begin{multline*}
%x^{16} + 8 x^{14} + 1432 x^{12} - 2752 x^{10} + 729648 x^8 - 5366016 x^6 + 204148608 x^4
%    - 1103153664 x^2 + 17167288576
%\end{multline*}
%is an unramified extension of $k$ with Galois group of 2-class 2 and with a degree 4 
%subextension over $k$ for which the 2-Sylow subgroup of its ideal class group is equal
%to $[2,16]$. 
is an unramified extension of $k$ with the 2-Sylow subgroup of its ideal class group 
equal to $[2,16]$. 
It follows that $G$ must be isomorphic to one of these two groups. 
Each one of these can be obtained by selecting $r \in \{ 0,1 \}$
and then considering the group generated by $\{x_i\}_{i=1}^8$ 
subject to the power commutator presentation
\begin{align*}
    [x_2,x_1] &= x_3 \hspace*{2.5cm}  &[x_5,x_2] &= x_8 \hspace*{2.5cm} &x_1^2 &= x_4 \\
    [x_3,x_1] &= x_5 &[x_5,x_3] &= x_8 &x_2^2 &= x_5 x_7 \\
    [x_3,x_2] &= x_6 &[x_5,x_4] &= x_8 &x_3^2 &= x_6 x_7 \\
    [x_4,x_2] &= x_5 x_6 x_7 x_8 &[x_7,x_1] &= x_8  &x_6^2 &= x_8\\
    [x_4,x_3] &= x_7 &[x_7,x_2] &= x_8  & &\\
    [x_5,x_1] &= x_7 &x_4^2 &= x_8^r & & 
\end{align*}   
% old::   \quad\makebox[2cm][l]{$\text{where}\:i \in \{0,1\}$}
Remark \ref{pcrem} about power commutator presentations applies here also.

Once again these groups are very similar. They both have order $2^8$ and 2-class 5
and by looking at their derived series we deduce
\begin{proposition} The field $k = \mathbb{Q}(\sqrt{-445})$ has finite $2$-class
tower of length 3, ie. $k = k_0 \subset k_1 \subset k_2 \subset k_3 = k^{nr,2}$.
We have $\text{Gal}\,(k_1/k_0) \cong [2,4]$, $\text{Gal}\,(k_2/k_1) \cong [2,2,4]$
and $\text{Gal}\,(k_3/k_2) \cong [2]$
\end{proposition}

\begin{remark}
We note that the finiteness of the towers in the propositions above is purely
a group theoretic result. Any field $k$ for which there exists a degree 8 subextension
$L$ of $k^{nr,2} / k$ such that the lattice of subfields and corresponding $2$-class
groups are the same as those in Figures \ref{fig2} and \ref{fig3} must also 
have a finite $2$-class tower. 
\end{remark}

\begin{remark}
\label{mult}
In each of the previous examples a single unramified degree 8 extension of $k$ was selected.
The number theoretic information provided by this extension was sufficient for us to be
able to show the finiteness of $\text{Gal}\,(k^{nr,2}/k)$ in each case. One obvious way
to extend the algorithm would be to incorporate the information from several such extensions
at the same time. In some sense we were already doing this (in an ad hoc fashion) when
we took the initial finite list of candidates in each of the previous examples 
and eliminated some of them using additional number theoretic information. More precisely
we now consider pairs of the form $(H,\{\{H^{(j)}_i\}_{i=1}^n\}_{j=1}^t)$ where $t$ is the number
of extensions under consideration. Such a  pair will be called a {\em representative} of $G_{(m)}$ 
if there exists a family of isomorphisms
$\psi_j : H \longrightarrow G/P_m(G)$ such that $\psi_j(H^{(j)}_i) = \overline{G_i}$
for each $i = 1,\ldots\,,n$ and $j = 1,\ldots\,,t$. As before given a list of pairs 
$\mathbf{L}_m$ containing a representative of $G_{(m)}$ we can use the p-group generation 
algorithm to compute a list $\mathbf{L}_{m+1}$ containing a representative of $G_{(m+1)}$. 
\end{remark}

\subsection{The fields $k = \mathbb{Q}(\sqrt{-1015})$ and $\mathbb{Q}(\sqrt{-1595})$}
In both cases $G/G' \cong [2,8]$ so we have $G/P_1(G) \cong [2,2]$.
The fields $\mathbb{Q}(\sqrt{-7},\sqrt{5},\sqrt{29})$ and 
$\mathbb{Q}(\sqrt{-11},\sqrt{5},\sqrt{29})$ are degree 4 subextensions
of $k^{nr,2}/k$ in each case. We can calculate the 2-Sylow subgroup of 
the ideal class group of each intermediate field, and construct a representative
of the pair $G_{(1)} = (G/P_1(G),\{\overline{G_i}\}_{i=1}^5)$ where the groups $G_i$ 
are the subgroups of $G$ corresponding to the intermediate fields. For both the given
fields we get the same representative of $G_{(1)}$.
Applying the method from Section 2 we are able to narrow down $G/P_2(G)$ to
one of 4 groups. Two of these can be identified as $[2,4]$ and $D_4$ (the 
dihedral group of order 8). These can be eliminated as possibilities since
it is easy to find examples of degree 8 subextensions of $k^{nr,2}/k$ with either of
these groups as the Galois group (see below). This implies that both groups must arise
as quotients of $G/P_2(G)$. It follows that $G/P_2(G)$ must be isomorphic to one
of the remaining two groups of order 16 which we will denote $H_3$ and $H_4$. The group 
$H_3$ is actually isomorphic to the group $H_3$ in the previous example.

\begin{figure}[h]
 \scalebox{.55}{\includegraphics{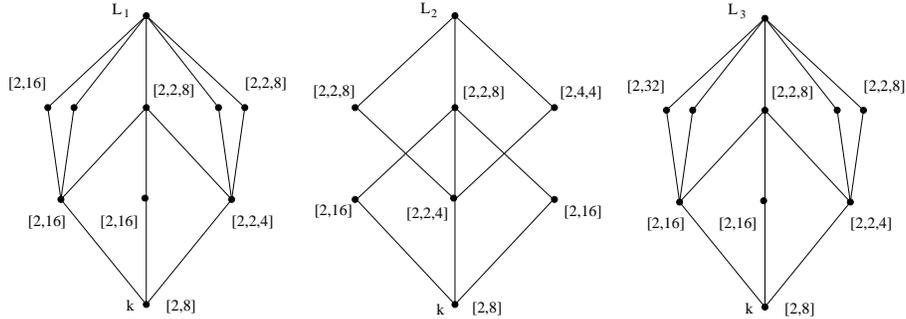}}
 \caption{Three subfield lattices and associated 2-class groups}\label{fig4}
\end{figure}

Now taking into account Remark~\ref{mult} we consider several unramified degree 8 extensions $L$ over $k$ instead
of only one. For $k = \mathbb{Q}(\sqrt{-1015})$ we consider the fields defined over
$\mathbb{Q}$ by the polynomials
\begin{eqnarray*}
 &L_1 :\: &x^{16} + 8302 x^{14} + 29865815 x^{12} + 60621449422 x^{10} + 75762817738769 x^8 \\
    & & +  59625975137422568 x^6  + 28858765154851072400 x^4 \\ 
    & & + 7861191091575524181248 x^2 + 924182332972720716353536 \\ 
 &L_2 :\: &x^{16} + 68 x^{14} - 26 x^{13} + 1922 x^{12} - 2316 x^{11} + 29806 x^{10} - 20958 x^9 + 
    335885 x^8 \\ 
    & & + 62002 x^7 + 1639268 x^6  + 2747082 x^5 + 6227217 x^4 + 
    7583004 x^3 + 7628823 x^2 \\ 
    & & + 4664142 x + 1486431  
\end{eqnarray*}
\begin{eqnarray*}
 &L_3 :\: &x^{16} + 208250 x^{14} + 6454080 x^{13} + 84986985877 x^{12} - 319881524440 x^{11} \\
   & & + 10504186175856042 x^{10} + 3217249977395280 x^9 
    + 2231207353583759168404 x^8 \\ 
    & & - 41652321975526297906680 x^7 + 162122207446267254901910082 x^6 \\ 
    & & -   4449276375660698756114160120 x^5 + 23169896204558457954443037721749 x^4 \\
    & & - 466821086268574071299245171753200 x^3  \\ 
   & & + 796726493601047682437367297422156178 x^2 \\
   & & - 37902470874562381569855570165787116760 x \\
   & & + 81595985378826852513556627342521162440521 \\
\end{eqnarray*}
%C35_29ck/kashdata --> list[10],[14],[16]
For $k =\mathbb{Q}(\sqrt{-1595})$ we consider the fields defined over
$\mathbb{Q}$ by the polynomials
\begin{eqnarray*}
 &L_1 :\: &x^{16} + 75 x^{14} + 3384 x^{12} + 85875 x^{10} + 1497421 x^8 + 16831500 x^6 \\ 
     & & + 129999744 x^4  + 564715200 x^2 + 1475789056 \\
 &L_2 :\: &x^{16} + 145 x^{14} + 4721 x^{12} - 336690 x^{10} + 2932126 x^8 + 22696270 x^6 \\
     & & + 32760881 x^4 - 90377775 x^2 + 43046721 \\
 &L_3 :\: &x^{16} - 2312 x^{14} + 2359542 x^{12} - 1183214812 x^{10} + 276742820433 x^8 \\
 	& & - 63131144780036 x^6 + 66951555767033248 x^4 - 13918403354887798784 x^2 \\ 
   & & + 793394183478882017536
\end{eqnarray*} 
%C55_29ck/kashdata --> list[18],[5],[16]
In each case we obtain three lattices of subfields 
together with the 2-Sylow subgroups of the corresponding
ideal class groups. These turn out to be the same for both $k = \mathbb{Q}(\sqrt{-1015})$ and
 $k =\mathbb{Q}(\sqrt{-1595})$ and are displayed in Figure \ref{fig4}. 
Note that $\text{Gal}\,(L_1/k) = D_4$, $\text{Gal}\,(L_2/k) = [2,4]$ and
$\text{Gal}\,(L_3/k) = D_4$.
A list $\mathbf{L}_2$ containing a representative of $G_{(2)}$ 
(see the modified definition in Remark~\ref{mult})
can be constructed by taking each group $H_r$ ($r = 3$ or $r = 4$) and adjoining 
all pairs of the form
$(H_r, \{\{f_j^{-1}(R^{(j)}_i)\}_{i\in I_j}\}_{j=1}^3)$ where $f_j$ runs through all the surjective homomorphisms
from $H_r$ onto $\text{Gal}\,(L_j/k)$ and $\{R^{(j)}_i\}_{i\in I_j}$ are the subgroups of 
$\text{Gal}\,(L_j/k)$ of index at most 4. Here $I_j$ is simply an indexing set 
for the subgroups, its
size being determined by the structure of $\text{Gal}\,(L_j/k)$. 
When the method from Section 2 is applied to $\mathbf{L}_2$ 
the sequence of lists generated terminates and a total of 2 candidates for $G$ are found.
Each one of these can be obtained by selecting $r \in \{ 0,1 \}$
and then considering the group generated by $\{x_i\}_{i=1}^9$ 
subject to the power commutator presentation
\begin{align*}
    [x_2,x_1] &=  x_3   \hspace*{2.5cm}  &[x_5,x_2] &= x_9 \hspace*{2.5cm} &x_2^2 &= x_5   x_7 \\
    [x_3,x_1] &=  x_5        &[x_5,x_3] &=  x_9  	&x_3^2 &= x_7   x_8  x_9^{1-r} \\
    [x_3,x_2] &=  x_8   x_9^{1-r}      &[x_5,x_4] &=  x_9  &x_4^2 &= x_6 \\
    [x_4,x_2] &=  x_5   x_7   x_8   x_9^r      &[x_7,x_1] &=  x_9 &x_6^2 &= x_8 \\
    [x_4,x_3] &=  x_7     &[x_7,x_2] &=  x_9  &x_8^2 &= x_9  \\
    [x_5,x_1] &=  x_7      &x_1^2 &= x_4  & &
\end{align*}   
Remark \ref{pcrem} about power commutator presentations applies here also.
Both of these groups have order $2^9$ and 2-class 5. By looking at their derived series we deduce
\begin{proposition} The fields $k = \mathbb{Q}(\sqrt{-1015})$ and
$k = \mathbb{Q}(\sqrt{-1595})$ have finite $2$-class
tower of length 3, ie. $k = k_0 \subset k_1 \subset k_2 \subset k_3 = k^{nr,2}$.
In both cases we have $\text{Gal}\,(k_1/k_0) \cong [2,8]$, $\text{Gal}\,(k_2/k_1) \cong [2,2,4]$
and $\text{Gal}\,(k_3/k_2) \cong [2]$
\end{proposition}

\end{document}